\documentclass[10pt]{article}
\usepackage{amsfonts}
\usepackage{amsmath}
\usepackage{amssymb}
\usepackage{amscd}

\newcommand{\beq}{\begin{equation}}
\newcommand{\eeq}{\end{equation}}
\newcommand{\ba}{\begin{array}}
\newcommand{\ea}{\end{array}}
\newcommand{\beqa}{\begin{eqnarray}}
\newcommand{\eeqa}{\end{eqnarray}}
\newcommand{\nn}{\nonumber \\}

\def\cala{{\cal A}}

\def\be{\begin{equation}}
\def\ee{\end{equation}}

\def\R{{\hat{\mathsf R}}}
\def\id{{\rm{id}}}

\newcounter{Theorem}\newcounter{Remark}
\newcounter{Definition}
\newcounter{Example}\newcounter{Exercise}

\newenvironment{Theorem}[1][\bf Theorem \arabic{Theorem}]{%

        \refstepcounter{Theorem}\noindent\textbf{#1.}${}$\hspace{1pt}${}$\it}{}

\newenvironment{Lemma}[1][\bf Lemma \arabic{Theorem}]{%
        
        \refstepcounter{Theorem}\noindent\textbf{#1.}${}$\hspace{1pt}${}$\it}{}
        
        \newenvironment{Definition}[1][\bf Definition~\arabic{Definition}]{%
        
        \refstepcounter{Definition}\noindent\textbf{#1.}${}$\hspace{1pt}${}$}{}

\begin{document}

\title{DRINFELD-JIMBO\\QUANTUM LIE ALGEBRA}

\author{
Oleg Ogievetsky
\footnote{On leave of absence from P.N. Lebedev Physical Institute, 
Theoretical Department, Leninsky prospekt 53, 119991 Moscow, Russia} \\[.6em]
\it Centre de Physique Th\'eorique \\ \it Luminy, 13288 Marseille, France \\ oleg@cpt.univ-mrs.fr 
\\[1.8em]
Todor Popov \\[.6em]
\it Institute for Nuclear Research and Nuclear Energy\\ 
\it Bulgarian Academy of Sciences \\ 
\it Sofia, BG-1784, Bulgaria \\
 tpopov@inrne.bas.bg
}

\date{}
\maketitle

\begin{abstract}
Quantum Lie algebras related to multi-parametric Drinfeld--Jimbo $R$-matrices
of type $GL(m\vert n)$ are classified.
%\keywords{Quantum groups; bicovariant  diffferential calculus; quantum Lie algebras.}
\end{abstract}

\section{Introduction}	

Quantum groups are examples of non-commutative manifolds having a rich differential geometry.
Woronowicz developed a general theory of the differential calculus on a quantum group, the so called bicovariant differential calculus.
Bicovariant bimodules are objects analogous to tensor bundles over Lie groups \cite{W}. The vector space dual to the space of left-invariant differential forms (i.e., the space of the left-invariant  vector fields) is endowed with a bilinear operation playing the role of the Lie bracket (for a friendly 
introduction to the subject see, e.g., \cite{AC,AS}). This vector space is a quantum analogue of a Lie algebra.
A quantum Lie algebra can be defined axiomatically as 
a triple $(V,\sigma, C)$ consisting of a vector space $V$, a braiding $\sigma: V\otimes V \rightarrow  V\otimes V$  
and a  ``quantum Lie bracket'' $C: V \otimes V \rightarrow V$ on it satisfying certain compatibility identities, see below. 

Given a braiding $\sigma$, it is natural to try to describe quantum Lie algebras compatible with 
$\sigma$. An important class of braidings arises as quantizations \cite{ESS,IOq} of classical $r$-matrices corresponding
to Belavin--Drinfeld triples \cite{BD}. The problem of a description of quantum Lie algebras compatible with the Cremmer--Gervais
$R$-matrix \cite{CG} has been addressed in our previous work \cite{OP2} with the help of a suitable rime Ansatz \cite{OP}. 
The Cremmer--Gervais $R$-matrix corresponds to a maximal Belavin--Drinfeld triple for 
the defining fundamental representation of the quantum group of type GL. In the defining representation this braiding is of Hecke type
which, in particular, implies that the BRST operator for a quantum Lie algebra with this braiding is finite \cite{IO2,IO}. 
In this note we give a complete list of quantum Lie algebras compatible with the Drinfeld-Jimbo
$R$-matrix \cite{D,J} (which corresponds to the empty Belavin--Drinfeld triple) in the defining fundamental representation on a vector super-space.

\section{Bicovariant calculus} 

We give here a short extract from the differential calculus on quantum groups to motivate the notion of the quantum Lie 
algebra and related objects used in the sequel.

Let $\cala$ be the Hopf algebra  of functions on a quantum group. We denote by $\Delta$, $\epsilon$ and $S$ the
coproduct, the counit and the antipode of $\cala$. We use the Sweedler notation 
for the coproduct
\[ \Delta(a) = a_{(1)} \otimes a_{(2)}   \ . \]
Woronowicz \cite{W} introduced the notion of the bicovariant bimodule $\Gamma$ over $\cala$.
It is an $\cala$-bimodule endowed with two coactions
\beq
\Delta_L: \Gamma \rightarrow \cala \otimes \Gamma\ \ ,\ \ \Delta_R: \Gamma \rightarrow  \Gamma \otimes \cala \nonumber  
\eeq
satisfying certain compatibility conditions. One supposes that $\Gamma$ is a free left (and right) $\cala$-module 
admitting a basis $\omega^i$ formed by {\it left-invariant forms}, that is,
\beq
\Delta_L(\omega^i)= 1 \otimes \omega^i \ .\eeq
Elements $r^i_j\in \cala$ are defined by 
\beq\Delta_R(\omega^i)=  \omega^j \otimes r^i_j \ .
\eeq
Denote the Hopf algebra dual to $\cala$ by $\cala'$; it has the coproduct $\Delta'$, counit $\epsilon'$ and
antipode $S'$. There exists a set of elements  $f^i_j\in \cala'$ relating the left and the right actions 
\beq \omega^i b = (f^i_j \ast b )\omega^j :=  b_{(1)} f^i_j( b_{(2)})\omega^j \ .\eeq
The consistency implies 
\beq\label{hopf}\Delta' f^i_j = f^i_k \otimes f^k_j \ ,\ \epsilon'(f^i_j)= \delta^i_j\ , 
\ S' (f^i_k) f^k_j = \delta^i_j=  f^i_k S'(f^k_j)\ ,\eeq
\beq\Delta r^i_j = r^k_j \otimes r^i_k \ ,\ \epsilon(r^i_j)= \delta^i_j\ , 
\ S(r^k_j ) r^i_k = \delta^i_j= r^k_j  S( r^i_k) \ .\eeq

The differential (on functions) in the Woronowicz calculus is the map $d:\cala\to\Gamma$, given by
\beq d a = (\chi _i \ast a) \omega^i  \qquad \forall\ a \in \cala \ ,\eeq
where the elements $\chi _i \in \cala'$ form a basis of the free (left) $\cala$-module of the {\it left-invariant vector fields}. 
  
The dual $\cala^{\prime}$ of a (finite dimensional) Hopf algebra $\cala$ is again a Hopf algebra. The Leibniz rule implies the coproduct
\beq\label{hopff} \Delta' \chi_i = \chi_j \otimes f^j_i + 1 \otimes \chi_i \ .\eeq
and  the duality induces  $ \epsilon' (\chi_i) = 0 $, hence $S'(\chi_j) = - \chi_i S'(f^i_j) $.

The elements $\chi_i$ and $f^i_j$ satisfy the quadratic-linear relations 
\beq\label{bcc}\begin{array}{lll}
&\chi_i \chi_j -   \chi_k \chi_l \sigma^{kl}_{ij} = \chi_k  C^k_{ij} \ ,
&  \sigma^{ab}_{kl} f^k_i f^l_j =f^a_k f^b_l \sigma^{kl}_{ij}  \ ,\\[1em]
 & \chi_k f^a_l \sigma^{kl}_{ij} +  f^a_l C^l_{ij} = C^a_{kl}f^k_i f^l_j +f^a_i \chi_j\ ,\qquad \qquad 
&\chi_i f^a_j =  f^a_k \chi_l  \sigma^{kl}_{ij}\ .  \end{array}\eeq
Here 
\beq\label{aut}  \sigma^{ij}_{kl}= f^i_l (r^j_k)\ \ \text{and}\ \ C^k_{ij} = \chi_j (r^k_i)  \ .\eeq 
The compatibility leads to the following relations between the {\it braiding} $\sigma$ and {\it structure constants} $C$:
\beq\label{braid}
\sigma^{ab}_{lm}\sigma^{mc}_{nk}\sigma^{ln}_{ij}= \sigma^{bc}_{lm} \sigma^{al}_{in}\sigma^{nm}_{jk}\qquad 
\mbox{(braid relation)}\ ,\eeq
\beq  C^b_{sk} C^s_{ij}   =   C^b_{is} C^s_{jk}  +  C_{sl}^b C_{ir}^s \sigma^{rl}_{jk}\qquad \mbox{(braided Jacobi identity)}\ ,\eeq
\beq \label{adj}  \sigma^{ab}_{sk} C^s_{ij} =   C^b_{sl} \sigma^{as}_{ir} \sigma^{rl}_{jk} \ , \qquad \eeq
\vskip -.5cm
\beq \label{adj0}  \sigma^{ab}_{sl} C^s_{ir} \sigma^{rl}_{jk} +  \sigma^{ab}_{i l} C^l_{jk}
=C^a_{rl} \sigma^{lb}_{sk} \sigma^{rs}_{ij}  +  C^b_{sk} \sigma^{as}_{ij}\ . \eeq

Denote by $\mathcal W$ the algebra with the generators $\chi_i$ and $f^i_j$ and the defining relations (\ref{bcc}). 
The formulas (\ref{hopf}) and (\ref{hopff}) equip $\mathcal W$ with a Hopf algebra structure. For 
further details see the original work \cite{W}.

The relations (\ref{braid})-(\ref{adj0}) can be conveniently written \cite{B,B2}
as a single braid relation. Let us make a convention that
the small  indices $i, j, \ldots ,k$ run over a set ${\cal I }=\{ 1,\dots,\dim V\}$ and the capital indices $I,J, \ldots, K $ run over the set 
${\cal I}_0:=0\cup {\cal I} $. 

\vskip .2cm
\begin{Lemma}
\label{lebe}  Let $\R^{IJ}_{KL}$ be a  matrix  whose non-vanishing components are
\beq\label{conv}\R^{ij}_{kl}= \sigma^{ij}_{kl} \ , \quad \R^{0j}_{kl} =
C^j_{kl} \ , \quad  \R^{0A}_{B0} = \delta^A_B  \ ,\quad   \R^{A0}_{0B} = \delta^A_B \ .\eeq 
Then the system (\ref{braid})-(\ref{adj0}) is equivalent to the braid relation
\beq \R_{12} \R_{23} \R_{12} = \R_{23} \R_{12} \R_{23}\ ,\eeq
where $\R_{12}:=\R\otimes {\rm{Id}}$ and $\R_{23}:={\rm{Id}}\otimes\R$.
Let $T^I_J$ be the matrix with elements $T^i_j = f^i_j, T^0_j=\chi_j, T^0_0=1, T_0^i=0$,
\beq  T^I_J=\left(\ba{cc}1 & \chi_j \\ 0& f^i_j\ea \right)\ .\eeq
Then the relations (\ref{bcc}) of the algebra $\cal W$ take the concise form
\beq\label{birt}  \R^{AB}_{KL} T_I^K T_J^L =T_K^A T_L^B  \R_{IJ}^{KL}\ .\eeq 
The Hopf structure of $\cal W$ (see  eqs. (\ref{hopf}) and (\ref{hopff})) reads
\beq\ba{ccc}
& \Delta' T^I_J = T^I_ K \otimes T^K_J\ , \qquad \epsilon'(T^I_J)=\delta^I_J \ ,&\\[1em] 
& S'(T^I_K)T^K_J = \delta^I_J = T^I_K S'(T^K_J)\ .&\ea\eeq
\end{Lemma}

The complexity of the relations (\ref{bcc}) of the algebra $\cal W$ is hidden into the single matrix $\R$
and the matrix of generators $T^I_J$. 

Recall that the (right) adjoint representation of a Hopf algebra ${\cal H}$ on itself is defined by $ad_x(y):=S( x_{(1)}) y x_{(2)}$.
The image of the space $V$ in $\cal W$ is stable with respect to the adjoint action and one has
$$ad_{f^i_j}(\chi_a)=\chi_b\sigma^{ib}_{aj}\ ,\ ad_{\chi_i}(\chi_a)=\chi_bC^b_{ai}\ .$$
In this representation the defining relations (\ref{bcc}) turn into the (numerical) relations (\ref{braid})-(\ref{adj0}). 

\section{Quantum Lie algebra}
The notion of a quantum Lie algebra formalizes the properties of the subalgebra of ``vector fields'', generated by $\chi_i$, in $\cal W$. 
We give the precise definitions. Let $V$ be the vector space with the basis $\{\chi _i\}_{i=1,\dots,\dim V}$.
The space $V$ is endowed with a braiding operator, that is, the operator $\sigma: V \otimes V \rightarrow V \otimes V$
which satisfies 
$$ \sigma_{12}\sigma_{23}\sigma_{12}=\sigma_{23}\sigma_{12}\sigma_{23}\ .$$
We assume that $\sigma$ is semi-simple and has an eigenvalue 1. Denote by $P_{(1)}$ the projector of $\sigma$ corresponding to the eigenvalue 1.

\vskip .2cm
\begin{Definition}\hspace{-.15cm}\hspace{.15cm} A  quantum Lie algebra is a triple $(V,\sigma ,C)$ 
where $C$ is the "bracket" 
$C: V \otimes V \rightarrow V$
such that the following conditions hold:

\vskip .1cm
i) braided symmetry: the bracket is in the kernel of the projector $P_{(1)}$
\beq \label{brsym} C P_{(1)}  =0\ ,\eeq
 
ii) braided Jacobi identity
\beq C (  C \otimes \id )  =    C ( \id \otimes C )+   C( C  \otimes \id) \sigma_{23} \ ,\eeq

iii) additional, linear in $C$, identities
$$\ba{rcl}   \sigma (C \otimes \id) &= &     (\id \otimes C) \sigma_{12} \sigma_{23}\ ,\\[.5em] 
\sigma (C\otimes \id) \sigma_{23}  +  \sigma (\id \otimes C)&= &  (C\otimes \id) \sigma_{23} \sigma_{12}
+ (\id \otimes C)\sigma_{12} \ .\ea$$
\end{Definition}

\vskip .2cm
\begin{Definition}\hspace{-.15cm}\hspace{.15cm} 
The universal enveloping algebra $U_{\sigma ,C}(V)$ of the quantum Lie al\-gebra $(V,\sigma ,C)$
is the associative algebra with the generators $\chi _i, i=1,\dots,\dim V$, and the defining relations 
\beq\label{qLie}\chi_i \chi_j -  \chi_k \chi_l \sigma_{ij}^{kl} =  \chi_k C_{ij}^k \ .\eeq
\end{Definition}

\noindent {\bf Notation.} We shall often write the compatibility conditions (\ref{adj}) and (\ref{adj0}) between $\sigma$ and $C$ in the form 
\beq\label{sub}E_{\sigma,C}(i,j,k;a,b):=  \sigma^{ab}_{sk} C^s_{ij}  -   C^b_{s \,l} \sigma^{as}_{ir} \sigma^{rl}_{jk}\ ,\eeq
\beq\label{sub2}F_{\sigma,C}(i,j,k;a,b):=\sigma^{ab}_{sl} C^s_{ir} \sigma^{rl}_{jk} +  \sigma^{ab}_{i l} C^l_{jk}
-C^a_{rl} \sigma^{lb}_{sk} \sigma^{rs}_{ij}  -  C^b_{sk} \sigma^{as}_{ij}\ .\eeq

\section{Ice condition}

The Boltzmann weights of the 6-vertex model in 2D statistical mechanics  are subject to a
restriction known as ``ice" condition. Namely, the ice condition for the entries of an R-matrix $\hat{R}$ can be different from zero only if the set of the upper and the set of the lower indices coincide
\beq \label{ice}\hat{R}^{ij}_{kl}\neq 0\qquad\Rightarrow\qquad\{ i,j\}\equiv\{ k,l\}\ .\eeq
An ice matrix $\hat{R}\in {\rm End}( V\otimes V)$ has the form
$$\hat{R}_{ij}^{kl} = a_{ij}\delta_i^l \delta_j^k + b_{ij}  \delta_i^k \delta_j^l\ . $$
If the set ${\cal I }$ of indices cannot be split into a disjoint union of two subsets ${\cal I }'$ and ${\cal I }''$ such that $b_{ij}=0=b_{ji}$ whenever $i\in {\cal I }'$ 
and $j\in {\cal I }''$ then we say that the matrix $\hat{R}$ is indecomposable. Recall that an operator $\hat{R}\in {\rm End}( V\otimes V)$
is called skew-invertible if there exists an operator $\Psi\in {\rm End}( V\otimes V)$ such that $\Psi^{iu}_{jv}\hat{R}^{vk}_{ul}=\delta^i_l\delta^k_j$.
The characteristic function $\theta_{i>j}$ is defined to be $1$ when $i>j$ and zero otherwise (similarly for $\theta_{i<j}$).

\vskip .2cm
\begin{Lemma} (See \cite{QS} for the proof.) Let $\hat{R}$ be an ice solution of the braid equation. Assume that $\hat{R}$ is invertible, skew-invertible and  indecomposable. Then, up to a reordering of the basis and rescaling, $\hat{R}$ is the standard multi-parametric Drinfeld-Jimbo R-matrix \cite{D,J} 
with $a_{ij}$  and $b_{ij}$ given by  
\beq\label{siceab}a_{ij}=(-1)^{\hat{i}\delta_{ij}} q^{(1-2\hat{i})\, \delta_{ij}}p_{ij}q^{\theta_{i<j}- \theta_{i>j}} \  , 
\quad  b_{ij}=(q-q^{-1}  )\theta_{i<j}\ . \eeq
Here $q\in \mathbb{C}^*$ and $\hat{i}\in \{ 0,1\}$; the parameters $p_{ij} $ satisfy $p_{ij} p_{ji}=1$ and $p_{ii}=1$.
\end{Lemma}

\section{Drinfeld-Jimbo quantum Lie algebra}
 
We shall say that $\hat{i}$ is the ``parity" of the basis vector $\chi_i$. Thus, $a_{ii}=q$ if  $\hat{i}=0$ and $a_{ii}=-q^{-1}$ if  $\hat{i}=1$.

\vskip .2cm
Recall that the braiding operator $\sigma$ of a quantum Lie algebra must have an eigenvalue 1. We thus have to rescale the standard R-matrix
whose eigenvalues are $q$ and $(-q^{-1})$; there are two possibilities:
$$\sigma = q^{-1}\hat{R} \qquad \mbox{or}  \qquad  \sigma = -  q\hat{R} \ .$$
The second possibility can be reduced to the first one. Namely, the replacement of $q$ by $(-\tilde{q}^{-1})$ in (\ref{siceab}) leads to the standard
R-matrix with the parities $\hat{i}'$, parameters $\tilde{q}=-q^{-1}$ and $\tilde{p}_{ij}$ where $\hat{i}'=1-\hat{i}$ and  
$\tilde{p}_{ij}=-p_{ij}\tilde{q}^{\, 2(\theta_{i>j}- \theta_{i<j})}$. Thus we have to investigate the quantum Lie algebras with the braiding operator 
$\sigma = q^{-1}\hat{R}$ or, explicitly, 
\beq\label{Siceab}\begin{array}{c}
\sigma_{ij}^{kl} = A_{ij}\delta_i^l \delta_j^k + B_{ij}  \delta_i^k \delta_j^l\ ,\\[1em]%\ 
A_{ij}=(-1)^{\hat{i}\delta_{ij}}q^{-2\hat{i}\, \delta_{ij}} p_{ij}q^{-2\theta_{i>j}}\ \ ,\ \ B_{ij}=1-q^{-2 \theta_{i<j}}\ .\end{array}\eeq
The R-matrix is called unitary if it squares to the identity operator. The matrix $\sigma$ is unitary iff $q^2=1$ and not semi-simple iff $q^2=-1$.

\vskip .2cm
\begin{Theorem}
\label{standQLA} Let $\sigma$ be a  standard Drinfeld-Jimbo R-matrix (\ref{Siceab}). Assume that
$\sigma$ is non-unitary and semi-simple, that is, $q^4\neq 1$. 
The non-trivial quantum Lie algebra with the braiding operator $\sigma$ exists only if the generator $\chi_1$ is even 
and $p_{1j}=1$ for $j>1$. It is then unique (up to a global rescaling of $C^k_{ij}$) and the structure constants $C^k_{ij}$ are given by 
\beq\label{pc}C^k_{ij}= c (\delta^1_i \delta^k_j - \delta^1_j \delta^k_i) \ .\eeq
\end{Theorem}
Explicitly, the  relations (\ref{qLie}) for this universal enveloping algebra $U_{\sigma ,C}(V)$ read 
\beq \label{qlice} \left\{\ba{rrcl}
 \chi_1 \chi_j - q^2 \chi_j \chi_1&=&q^2c  \chi_j \qquad&  \text{if}\ \ \ 1<j  \ , \nn
 \chi_i \chi_j - p_{ij} q^2 \chi_j \chi_i &=&0  \qquad& \text{if}\ \ \ 1< i<j  \ , \nn
\chi_i^2&=&0 \qquad& \text{if}\ \ \ \hat{i}=1\ .
 \ea \right.  
 \eeq

We start with the following lemma.

\vskip .2cm
\begin{Lemma}
\label{smallem} For the operator $\sigma$ given by (\ref{Siceab}) with $q^4\neq 1$, the relations (\ref{brsym}), (\ref{adj}) and (\ref{adj0})
imply
\beq C^k_{ji} = -  p_{ji} C^k_{ij} \qquad \mbox{and} \qquad C^k_{jj}=0 \ .\eeq
\end{Lemma}
{\bf Proof.} The relation (\ref{brsym}) is equivalent to 
$C^k_{ab} (\sigma^{ab}_{ij} + q^{-2}\delta^{a}_{i}\delta^{b}_{j}) =0$ for the
Hecke matrix $\sigma$ with the eigenvalues 1 and $(-q^{-2})$. For $i\neq j$, this immediately yields
$$  C^k_{ji} = -  p_{ji} C^k_{ij}\ \ , \ \ i \neq j \ .$$
{}For $i=j$ and $\hat{j}=0$ the relation (\ref{brsym}) reads $(1+q^{-2}) C^k_{jj} = 0$ so, by the semi-simplicity of $\sigma$, we obtain 
\beq\label{jj0} \qquad C^k_{jj}=0\ \ ,\ \  \hat{j}=0\ .\eeq
{}For $i=j$ and $\hat{j}=1$ the relation (\ref{brsym})
does not impose any constraint on $C^k_{jj}$
and we have to use eqs. (\ref{adj0}) and (\ref{adj}). 
Consider the equation $E_{\sigma,C}(j,j;j,j,j)=0$ with $\hat{j}=1$. The summation indices get fixed by the ice condition and we obtain
$$\sigma^{jj}_{jj} C^j_{jj}  =  C^j_{jj} \sigma^{jj}_{jj} \sigma^{jj}_{jj}\, \mbox{(no summation)}\Rightarrow 
 C^j_{jj} \left( 1- (-1)^{\hat{j}} q^{-2 \hat{j}}\right) =0  \, \Rightarrow  C^j_{jj} = 0 $$
by the semi-simplicity condition.
It is left to show that $C^k_{jj}$ for $\hat{j}=1$ and $j\neq k $. Choose  the following equations from the system (\ref{adj})-(\ref{adj0}):
\beqa E_{\sigma,C}(j,j,j;j,k)&=&0  \ ,\\   E_{\sigma,C}(j,j,k;k,k) &=&0   \ ,\\ 
F_{\sigma,C}(j,j,j;j,k)&=&0 \ .\eeqa 
Again the ice condition ``freezes" all summations and we obtain
\beq\ba{rll}  (\sigma^{jk}_{kj} - \sigma^{jj}_{jj}\sigma^{jj}_{jj})C^k_{jj} &=& 0\ ,    \\[.6em]
(\sigma^{kk}_{kk}-\sigma^{kj}_{jk} \sigma^{kj}_{jk}) C^k_{jj}&=&0\ ,    \\[.6em]
(\sigma^{jj}_{jj} \sigma^{jk}_{kj} + \sigma^{jk}_{jk} - \sigma^{jj}_{jj})C^k_{jj}&=& 0\ . 
 \ea  \qquad  \qquad \mbox{(no summation!)} \eeq 
Substituting the values of the matrix elements of  $\sigma$ we obtain (recall that $\hat{j}=1$) 
 \beq\label{koejjk} \ba{rll}  \left[\,  p_{kj}q^{- 2 \theta_{k>j}} -q^{-4 }\, \right]\, C^k_{jj} &=& 0\ ,\\[.6em]
\left [\,  (-1)^{\hat{k}}q^{-2 \hat{k}}- p^2_{jk} q^{-4 \theta_{j>k}}\, \right]\,  C^k_{jj} &=& 0\ , \\[.6em]
 \left[\,  -q^{-2}p_{kj} q^{ - 2\theta_{k>j}} +1- q^{-2 \theta_{j<k}}+ q^{-2}\, \right]\, C^k_{jj} & =&0\ . \ea \eeq
If $C^k_{jj}\neq 0$ then each square bracket  in (\ref{koejjk}) vanishes. It is straightforward to see that this contradicts to the 
restriction $q^4\neq 1$. \hfill $\Box$

\vskip .2cm
\noindent 
{\bf Proof of Theorem \ref{standQLA}.} We first prove that
$$C^k_{ij}=0\ \ , \quad i \neq j \neq k \neq i \ .$$ 
Consider the subsystem 
\beq\label{suyst}\ba{cccc}
E_{\sigma,C}(i,j,k;k,k)=0 &\Leftrightarrow&  A_{kk} C^k_{ij} =A_{jk} A_{ik} C^k_{ij}  ,& \\[.4em]
E_{\sigma,C}(i,j,j;j,k)=0&\Leftrightarrow& A_{kj} C^k_{ij} = A_{jj}A_{ij} C^k_{ij},& \mbox{(no summation!)} \\[.4em]
E_{\sigma,C}(i,j,i;i,k)=0&\Leftrightarrow& A_{ki} C^k_{ij} =A_{ji}A_{ii} C^k_{ij} &
\ea 
\eeq
with three different indices $i \neq j \neq k \neq i $. Since $A_{ij}A_{ji}=q^{-2}$ the system (\ref{suyst}) has a non-zero solution $C^k_{ij}$ iff 
$$  A_{ii} A_{jj} A_{kk} = q^{-2}\ \ \text{or}\ \ (-1)^{\hat{i}+\hat{j}+\hat{k}}q^{-2(\hat{i}+\hat{j}+\hat{k})}=q^{-2}\ .$$
This equation may have a solution only if $q^4=-1$ and $\hat{i}=\hat{j}=\hat{k}=1$. But then the equation $F_{\sigma,C}(i,i,j;i,k)=0$ implies
that $C^k_{ij}=0$. 
\vskip .2cm

Next we consider the equation $E_{\sigma,C}(k,i,j;i,k)=0$ with $i \neq j \neq k \neq i $; by the ice condition, this equation reduces to
$0= \sigma^{ij}_{ij} \sigma^{ik}_{ki} C^k_{k \,j} $. Since $ \sigma^{ik}_{ki} \neq 0$ and $\sigma^{ij}_{ij}\neq 0$ iff $i<j$ we conclude that 
$C^k_{kj} = 0$ iff $i<j$. If $j>1$ we can always find $i\colon i<j$. 
Thus the possibly non-zero structure constants $C^k_{kj}$  
are only $C^k_{k1}$; the structure constants $C^k_{1k}$ may also be different from 0, $C^k_{1k}= - p_{1k}C^k_{k1}$ by Lemma \ref{smallem}.
As we have seen, all other structure constants vanish.

The  constants $C^k_{1k}$ are subject to further constraints resulting from the equation $E_{\sigma,C}(1,j,k;j,k)=0$ with $1<j<k$,
\beq\label{cjk} \sigma^{jk}_{jk} C^j_{1j}  =     C^k_{1 k} \sigma^{j1}_{1j} \sigma^{jk}_{jk} \ 
\Rightarrow\ C^j_{1j}= p_{1j} C^k_{1 k}\ \ \mbox{for}\ \ j,k\ \ \mbox{such that}\ \ 1<j<k\eeq
(no summation). The relation $E_{\sigma,C}(i,j,j;j,j)=0$ reads
$$\sigma^{jj}_{jj} C^j_{ij}  =   C^j_{ij} \sigma^{ji}_{ij} \sigma^{jj}_{jj}\qquad \mbox{(no summation)}\ ,$$
which implies for  $i=1$ that
$$(1- p_{1j})C^j_{1j} =0\ .$$
Thus if $p_{1j} \neq 1$ for some $j$ then $C^j_{1j} =0$ and, by (\ref{cjk}), $C^k_{1k}=0$ for all $k$ and there is no non-trivial quantum Lie algebra. 
Therefore $p_{1j} =1$ and then, by (\ref{cjk}), the constants $C^k_{1k}$ for all $k\neq 1$ are equal, $C^k_{1k}=c$.
Now consider the equation $E_{\sigma,C}(1,j,1;1,j)=0$, 
$$C_{1j}^j \sigma^{1j}_{j1}(1 - \sigma^{11}_{11} ) = 0 \qquad \mbox{(no summation)}\ ,$$ 
So, if $\chi_1$ is odd then all structure constants vanish. This is our final result: 
\beq\label{fin}p_{1j} =1\ \  ,\ \  C^k_{ij} = c (\delta^k_j \delta^1_i -  \delta^k_i \delta^1_j)\ \ ,\ \ \chi_1\ \ \text{is even}\ .\eeq
It is straightforward to check that (\ref{fin}) defines the quantum Lie algebra. The proof is finished.\hfill $\Box$

\vskip .2cm
The braid relation is stable under three canonical operations (and their compositions): $\sigma\mapsto\sigma^T$ (transposition), 
$\sigma\mapsto\sigma_{21}:=P\sigma P$, where $P\in \text{End}(V\otimes V)$ is the flip, $P(u\otimes v)=v\otimes u$, and $\sigma\mapsto\sigma^{-1}$. 
However for the standard R-matrix $\sigma$, 
the matrices $\sigma^T$, $\sigma_{21}$ and $\sigma^{-1}$ are again standard (modulo a base change and a redefinition of parameters) and 
we do not need to consider them separately. We conclude that Theorem \ref{standQLA} gives all quantum Lie algebras compatible with the
standard Drinfeld-Jimbo $R$-matrices. 

\subsection*{Acknowledgement}
T.P. thanks for the warm hospitality the organizers of the SEENET-MTP Workshop in memory of Julius Wess. 
The work  was  supported by the Bulgarian NSF grant DO 02-257 and
the French-Bulgarian collaboration project ``Rila-4" N112.


\begin{thebibliography}{00}

\bibitem{AC} P. Aschieri and L. Castellani, {\em An introduction to noncommutative differential geometry on quantum groups}, Int. J. Mod. Phys. A{\bf 8} (1993) 1667--1706.

\bibitem{AS}P. Aschieri and P. Schupp, {\em Vector Fields on Quantum Groups},	Int. J. Mod. Phys. A{\bf 11} (1996) 1077--1100.

\bibitem{BD} A. A. Belavin and V. G. Drinfeld, {\em Solutions of the classical Yang--Baxter equation for simple Lie
algebras}, Funct. Anal. Appl. {\bf 16} (1982) 159--180.

\bibitem{B} D. Bernard, {\em Quantum Lie algebras and differential calculus on quantum groups}, 
Prog. Theor. Phys. Suppl. {\bf 102} (1990) 49-66.

\bibitem{B2} D. Bernard, {\em A remark on quasitriangular quantum Lie algebras}, Phys. Lett. {\bf B260} (1991) 389--393. 

\bibitem{CG} E. Cremmer and J.-L. Gervais, {\em The quantum group structure
associated with non-linearly extended Virasoro algebras};
Comm. Math. Phys. {\bf 134 } (1990) 619--632.

\bibitem{D} V. Drinfeld, {\em Quantum groups}, Proc. Int. Congr. Math. {\bf 1} (1986) 798--820. Amer. Math. Soc., Providence, RI, 1987.

\bibitem{ESS} P. Etingof, T. Schedler and O. Schiffmann, {\em Explicit quantization of dynamical r-matrices forfinite dimensional semisimple Lie algebras}, J. Amer. Math. Soc. {\bf 13} (2000) 595--609.arXiv: math.QA/9912009

\bibitem{IO2} A. Isaev and O. Ogievetsky, {\em BRST operator for quantum Lie algebras: explicit formula}, Int. J. Mod. Phys. A{\bf 19} (2004) 240--247. 

\bibitem{IO} A. P. Isaev and O. V. Ogievetsky, {\em BRST operator for quantum Lie algebras and differential calculus on 
quantum groups}, Theor. and Math. Phys. {\bf 129} (2001) 1558--1572.  arXiv: math/0106206

\bibitem{IOq} A. Isaev and O. Ogievetsky, {\em On quantization of r-matrices for Belavin--Drinfeld triples}, Phys. Atomic Nuclei {\bf 64} (2001) 2126--2130. 
arXiv: math.QA/0010190

\bibitem{J} M. Jimbo, {\em A $q$-difference analogue of $U_q(g)$ and the Yang--Baxter equation}, Lett. Math. Phys. {\bf 10} (1985) 63--69.
  
\bibitem{QS}  O. Ogievetsky, {\em Uses of quantum spaces}, Contemp. Math. {\bf 294} (2002) 161--232.

\bibitem{OP} O. Ogievetsky and T. Popov, {\em R-matrices in rime}, Adv. Theor. Math. Phys. {\bf 14} (2010) 439--506. arXiv: 0704.1947 [math.QA]

\bibitem{OP2} O. Ogievetsky and T. Popov, {\em Cremmer--Gervais quantum Lie algebra}, Fortsch. Phys.
{\bf 57} (2009) 654--658. arXiv: 0905.0882 [math-ph]

\bibitem{W} S. L. Woronowicz, {\em Differential calculus on quantum matrix pseudogroups
(quantum groups)}, Comm. Math. Phys. {\bf 122} (1989) 125--170.
\end{thebibliography}
\end{document}